\documentclass[11pt]{article} 
\usepackage{latexsym} 
\usepackage{amsmath}
\usepackage{amssymb} 
\usepackage{amsfonts}
\usepackage{amscd} 
\usepackage{graphicx}
\usepackage{layout}
\begin{document} 
%%%%%%%%%%%%%%%%%%%%%%%%%%%%%%%%% 
%%%%%%%%%%%%%%%%%%%%%%%%%%%%%%%%%
\newtheorem{Th}{Theorem}[section]
\newtheorem{Cor}{Corollary}[section]
\newtheorem{Prop}{Proposition}[section]
\newtheorem{Lem}{Lemma}[section]
\newtheorem{Def}{Definition}[section]
\newtheorem{Rem}{Remark}[section]
\newtheorem{Ex}{Example}[section]
\newtheorem{stw}{Proposition}[section]

%Definitions of bet ent

\newcommand{\bet}{\begin{Th}}
\newcommand{\ent}{\stepcounter{Cor}
   \stepcounter{Prop}\stepcounter{Lem}\stepcounter{Def}
   \stepcounter{Rem}\stepcounter{Ex}\end{Th}}

%Definitions of bec enc bep enp bel enl
%bef enf ber enr bee ene

\newcommand{\bec}{\begin{Cor}}
\newcommand{\enc}{\stepcounter{Th}
   \stepcounter{Prop}\stepcounter{Lem}\stepcounter{Def}
   \stepcounter{Rem}\stepcounter{Ex}\end{Cor}}
\newcommand{\bep}{\begin{Prop}}
\newcommand{\enp}{\stepcounter{Th}
   \stepcounter{Cor}\stepcounter{Lem}\stepcounter{Def}
   \stepcounter{Rem}\stepcounter{Ex}\end{Prop}}
\newcommand{\bel}{\begin{Lem}}
\newcommand{\enl}{\stepcounter{Th}
   \stepcounter{Cor}\stepcounter{Prop}\stepcounter{Def}
   \stepcounter{Rem}\stepcounter{Ex}\end{Lem}}
\newcommand{\bef}{\begin{Def}}
\newcommand{\enf}{\stepcounter{Th}
   \stepcounter{Cor}\stepcounter{Prop}\stepcounter{Lem}
   \stepcounter{Rem}\stepcounter{Ex}\end{Def}}
\newcommand{\ber}{\begin{Rem}}
\newcommand{\enr}{
   %\stepcounter{Rem} 
   \stepcounter{Th}\stepcounter{Cor}\stepcounter{Prop}
   \stepcounter{Lem}\stepcounter{Def}\stepcounter{Ex}\end{Rem}}
\newcommand{\bee}{\begin{Ex}}
\newcommand{\ene}{
 %\stepcounter{Ex}
   \stepcounter{Th}\stepcounter{Cor}\stepcounter{Prop}
   \stepcounter{Lem}\stepcounter{Def}\stepcounter{Rem}\end{Ex}}
\newcommand{\Proof}{\noindent{\it Proof\,}:\ }
\newcommand{\beP}{\Proof}
\newcommand{\enP}{\hfill $\Box$ \par\vspace{5truemm}}
%\newcommand{\enP}{\QED}

%%%%%%%%%%%%%%%%%%%%%%%%%%%%%%%%%%%%%%%%%
%Beginning of Local Definition
%Local definitions
\newcommand{\EE}{\mathbf{E}}
\newcommand{\QQ}{\mathbf{Q}}
\newcommand{\R}{\mathbf{R}}
\newcommand{\C}{\mathbf{C}}
\newcommand{\ZZ}{\mathbf{Z}}
\newcommand{\KK}{\mathbf{K}}
\newcommand{\NN}{\mathbf{N}}
\newcommand{\PP}{\mathbf{P}}
\newcommand{\HH}{\mathbf{H}}
\newcommand{\uuu}{\boldsymbol{u}}
\newcommand{\xxx}{\boldsymbol{x}}
\newcommand{\aaa}{\boldsymbol{a}}
\newcommand{\bbb}{\boldsymbol{b}}
\newcommand{\AAA}{\mathbf{A}}
\newcommand{\BBB}{\mathbf{B}}
\newcommand{\ccc}{\boldsymbol{c}}
\newcommand{\iii}{\boldsymbol{i}}
\newcommand{\jjj}{\boldsymbol{j}}
\newcommand{\kkk}{\boldsymbol{k}}
\newcommand{\rrr}{\boldsymbol{r}}
\newcommand{\FFF}{\boldsymbol{F}}
\newcommand{\yyy}{\boldsymbol{y}}
\newcommand{\ppp}{\boldsymbol{p}}
\newcommand{\qqq}{\boldsymbol{q}}
\newcommand{\nnn}{\boldsymbol{n}}
\newcommand{\vvv}{\boldsymbol{v}}
\newcommand{\eee}{\boldsymbol{e}}
\newcommand{\fff}{\boldsymbol{f}}
\newcommand{\www}{\boldsymbol{w}}
\newcommand{\0}{\boldsymbol{0}}
\newcommand{\lon}{\longrightarrow}
\newcommand{\ga}{\gamma}
\newcommand{\pa}{\partial}
\newcommand{\QED}{\hfill $\Box$}
\newcommand{\id}{{\mbox {\rm id}}}
\newcommand{\Ker}{{\mbox {\rm Ker}}}
\newcommand{\grad}{{\mbox {\rm grad}}}
\newcommand{\ind}{{\mbox {\rm ind}}}
\newcommand{\rot}{{\mbox {\rm rot}}}
\newcommand{\diver}{{\mbox {\rm div}}}
\newcommand{\Gr}{{\mbox {\rm Gr}}}
\newcommand{\LG}{{\mbox {\rm LG}}}
\newcommand{\Diff}{{\mbox {\rm Diff}}}
\newcommand{\Symp}{{\mbox {\rm Symp}}}
\newcommand{\Ct}{{\mbox {\rm Ct}}}
\newcommand{\Uns}{{\mbox {\rm Uns}}}
\newcommand{\rank}{{\mbox {\rm rank}}}
\newcommand{\sign}{{\mbox {\rm sign}}}
\newcommand{\Spin}{{\mbox {\rm Spin}}}
\newcommand{\Sp}{{\mbox {\rm sp}}}
\newcommand{\Int}{{\mbox {\rm Int}}}
\newcommand{\Hom}{{\mbox {\rm Hom}}}
\newcommand{\Tan}{{\mbox {\rm Tan}}}
\newcommand{\codim}{{\mbox {\rm codim}}}
\newcommand{\ord}{{\mbox {\rm ord}}}
\newcommand{\Iso}{{\mbox {\rm Iso}}}
\newcommand{\corank}{{\mbox {\rm corank}}}
\def\mod{{\mbox {\rm mod}}}
\newcommand{\pt}{{\mbox {\rm pt}}}
\newcommand{\qed}{\hfill $\Box$ \par}
\newcommand{\spe}{\vspace{0.4truecm}}
\newcommand{\ad}{{\mbox{\rm ad}}}

\newcommand{\dint}[2]{{\displaystyle\int}_{{\hspace{-1.9truemm}}{#1}}^{#2}}

%\newcommand{\dfrac}[2]{\frac{\displaystyle{#1}}{\displaystyle{#2}}}
%%%%%%%%%%%%%
%枠で囲むマクロ
%
\newenvironment{FRAME}{\begin{trivlist}\item[]
	\hrule
	\hbox to \linewidth\bgroup
		\advance\linewidth by -10pt
		\hsize=\linewidth
		\vrule\hfill
		\vbox\bgroup
			\vskip5pt
			\def\thempfootnote{\arabic{mpfootnote}}
			\begin{minipage}{\linewidth}}{%
			\end{minipage}\vskip5pt
		\egroup\hfill\vrule
	\egroup\hrule
	\end{trivlist}}
%%%%%%%%%%%%%%%%%%%%%%%%%%%%%%%%%%%%%%%%%%%%%%%%%
%End of local definitions
%%%%%%%%%%%%%%%%%%%%%%%%%%%%%%%%%%%%%%%%%
%End of Local definitions 
%%%%%%%%%%%%%%%%%%%%%%%%%%%%%%%%%

\title{
Singularities 
of tangent surfaces to generic space curves
} 

\author{G. Ishikawa
%\thanks{This work was supported by KAKENHI No.22340030 and No.23654058.}, 
and T. Yamashita}

%\renewcommand{\thefootnote}{\fnsymbol{footnote}}
%\footnotetext{Key words: 
%}
%\footnotetext{2000 {\it Mathematics Subject Classification}:  
%Primary 53C17; Secondly 58A30, 57R45, 93B05. }

\date{ }

\maketitle

%
%\begin{abstract} 
%
%\end{abstract}
%

\section{Introduction}
\label{Introduction}

The tangent lines to a space curve form a ruled surface, which is called the {\it tangent surface} 
or the 
{\it tangent developable} or the {\it tangential variety} to the space curve 
(see for instance \cite{Kreyszig}\cite{Cayley}\cite{Ishikawa4}\cite{Lawrence}). 
Tangent surfaces appear in various geometric problems and applications naturally, 
providing several important examples of non-isolated singularities in applications of geometry
(see for instance \cite{AG}\cite{CI}\cite{FP}\cite{IMT0}\cite{IKY}\cite{INS}\cite{Nuno Ballesteros}\cite{NS}). 

It is known, in the three dimensional Euclidean space ${E}^3$, 
that the tangent surface to a generic space curve $\gamma : I \to {E}^3$ 
is locally diffeomorphic to the {\it cuspidal edge}
or to the {\it folded umbrella} (also called, {\it cuspidal cross cap}), as is found by Cayley and Cleave \cite{Cleave}. 
Cuspidal edge singularities appear along ordinary points where 
$\gamma', \gamma'', \gamma'''$ are linearly independent, while 
the folded umbrella appears at an isolated point of zero torsion where 
$\gamma', \gamma'', \gamma'''$ are linearly dependent but $\gamma', \gamma'', \gamma''''$ 
are linearly independent (see \cite{BG}\cite{Porteous}). 

The classification is generalized to more degenerate cases 
by Mond \cite{Mond1}\cite{Mond3} and Shcherbak \cite{Shcherbak1}\cite{Shcherbak2}\cite{Arnol'd3}. 
See also \cite{Ishikawa3}\cite{Ishikawa4}. 
The classifications were performed mainly in locally projectively flat cases so far. 
However more general cases, namely, not necessarily projectively flat cases have never been treated as far as the authors recognize. 

\

In this paper we give the complete solution to 
the local diffeomorphism classification problem of generic singularities which appear in 
tangent surfaces, in as wider situations as possible. 
Naturally we interpret {\it geodesics} as \lq\lq lines" whenever a (semi-)Riemannian metric, 
or, more generally, an affine connection $\nabla$ is given in an ambient space of arbitrary dimension. Then, given 
an immersed curve, we define {\it $\nabla$-tangent surface}  
by the surface ruled by tangent geodesics to the curve. 
Then the main theorems in this paper are as follows: 

\bet
\label{genericity-theorem} {\rm (Genericity: Singularities of tangent surfaces to generic curves)}
Let $\nabla$ be any affine connection on a manifold $M$ of dimension $m \geq 3$. 
The singularities of the $\nabla$-tangent surface 
to a generic curve in $M$ on a neighborhood of the curve 
are only the cuspidal edges and the folded umbrellas if $m = 3$, 
and the embedded cuspidal edges if $m \geq 4$. 
\ent

%\bet
%\label{genericity-theorem2} {\rm (Genericity 2: Singularities of tangent surfaces to generic directed curves) }
%Let $\nabla$ be any affine connection on a manifold $M$ of dimension $m \geq 3$. 
%The singularities of the $\nabla$-tangent surface to a generic directed 
%curve in $M$ on a neighborhood of the curve 
%are only the cuspidal edges, the folded umbrellas and the swallowtails if $m = 3$,  
%and the embedded cuspidal edges and open swallowtails if $m \geq 4$. 
%\ent

The genericity is exactly given (see Propositions \ref{genericity1}) using 
Whitney $C^\infty$ topology on appropriate space of curves. 

\bet
\label{characterization-theorem}
{\rm (Characterization)} 
Let $\nabla$ be a torsion free affine connection on a manifold $M$. 
Let $\gamma : I \to M$ be a $C^\infty$ curve from an open interval $I$. 

{\rm (1)}  Let $\dim(M) = 3$. 
If $(\nabla\gamma)(t_0), (\nabla^2\gamma)(t_0), (\nabla^3\gamma)(t_0)$ are linearly independent, 
then the $\nabla$-tangent surface 
$\nabla{\mbox{\rm -}}\Tan(\gamma)$ is locally diffeomorphic to the cuspidal edge at $(t_0, 0) \in I \times \R$. 
If $(\nabla\gamma)(t_0), (\nabla^2\gamma)(t_0), (\nabla^3\gamma)(t_0)$ are linearly dependent, 
and $(\nabla\gamma)(t_0), (\nabla^2\gamma)(t_0), (\nabla^4\gamma)(t_0)$ are linearly independent, then
$\nabla{\mbox{\rm -}}\Tan(\gamma)$ is 
locally diffeomorphic to the folded umbrella at $(t_0, 0) \in I \times \R$. 

{\rm (2)} Let $\dim(M) \geq 4$. 
If $(\nabla\gamma)(t_0), (\nabla^2\gamma)(t_0), (\nabla^3\gamma)(t_0)$ are linearly independent, 
then the $\nabla$-tangent surface $\nabla{\mbox{\rm -}}\Tan(\gamma)$ is locally diffeomorphic to the embedded
cuspidal edge at $(t_0, 0) \in I \times \R$. 
\ent

\

A map-germ $f : (\R^2, p) \to M$ is locally {\it diffeomorphic} at $p$ to another map-germ 
$g : (\R^2, p') \to M'$
if there exist diffeomorphism-germs 
$\sigma : (\R^2, p) \to (\R^2, p')$ and $\tau : (M, f(p)) \to (M', g(p'))$ such that 
$\tau\circ f = g\circ \sigma : (\R^2, p) \to (M', g(p'))$. 

The {\it cuspidal edge} is defined by the map-germ $(\R^2, 0) \to (\R^m, 0)$, $m \geq 3$, 
$$
(t, s) \mapsto (t + s, \ t^2 + 2st, \ t^3 + 3st^2, \ 0, \ \dots, \ 0), 
$$
which is diffeomorphic to $(u, w) \mapsto (u, w^2, w^3, 0, \dots, 0)$. 
The cuspidal edge singularities are originally defined only in the three dimensional space. 
Here we are generalizing the notion of the cuspidal edge in higher dimensional space. In Theorem \ref{characterization-theorem} (2), we emphasize it by writing \lq\lq embedded" cuspidal edge. In what follows, we call it just 
cuspidal edge for simplicity even in the case $m \geq 4$. 
The {\it folded umbrella} (or the {\it cuspidal cross cap}) is defined by the map-germ 
$(\R^2, 0) \to (\R^3, 0)$, 
$$
(t, s) \mapsto (t + s, \ t^2 + 2st, \ t^4 + 4st^3), 
$$
which is diffeomorphic to $(u, t) \mapsto (u, t^2 + ut, t^4 + {\textstyle \frac{2}{3}}ut)$. 

\

In \cite{INS}, Izumiya, Nagai, Saji introduced and studied the class 
\lq\lq E-flat" great circular surfaces in the standard three sphere $S^3$ 
in detail, which contains the class of tangent surfaces to curves in $S^3$. 
The generic classifications given there (Theorems 1.2, 1.3 of \cite{INS}) in the sphere geometry 
become different from ours, 
because of the differences of topology and mappings spaces defining the genericity. 

Because we treat singularities in a general ambient space, we need the intrinsic 
characterizations of singularities found in \cite{FSUY}. 
The characterization of folded umbrellas is applied to Lorenz-Minkowski geometry in \cite{FSUY}. 
In this paper we apply to non-flat projective geometry the characterizations and their some generalization 
via the notion of openings introduced by the first author (\cite{Ishikawa4}, see also \cite{Ishikawa2}). 

\

In \S \ref{Affine connection and geodesic}, we recall on affine connections and related notions. 
In \S \ref{Tangent surface and frontal}, we define the tangent surface to an immersed curve and 
show that the tangent surface is a frontal under certain general conditions (Lemma \ref{non-degenerate-frontal}). 
We recall the criteria of singularities in \S \ref{Cuspidal edge and folded umbrella}. 
After a preliminary calculations in \S \ref{Characteristic vector field}, 
we show Theorem \ref{characterization-theorem} 
in \S \ref{Euclidean case}, using the criteria of singularities for the Euclidean case and 
in \S \ref{Proof of the characterization theorem} in general. 
In \S \ref{Proof of the genericity theorem 1} we show Theorem \ref{genericity-theorem}. 
Apart from main theorems, but related to them, we 
give an observation on the singularities of tangent surfaces to torsionless curves in 
\S \ref{Tangent surfaces to torsionless curves and fold singularities}. In flat case the 
tangent surface to a torsionless curve necessarily has the \lq\lq fold" singularity. However, in non-flat case 
we have an example where $(2, 5)$-cuspidal edge singularity appears 
on the tangent surface of some torsionless curve. 

\

In this paper all manifolds and mappings are assumed to be of class $C^\infty$ unless otherwise stated. 

\

This paper is derived as the main extract of the results in an unpublished paper \cite{IY}. 

\section{Affine connection and geodesic}
\label{Affine connection and geodesic}

Let $M$ be an $m$-dimensional manifold with an affine connection 
$\nabla$ (see \cite{Helgason}\cite{KN}). For any vector fields $X, Y$ on $M$, a 
vector field $\nabla_XY$ on $M$, which is called the {\it covariant derivative} of $Y$ by $X$, is assigned such that 
$$
\nabla_{hX + kY}Z = h\nabla_XZ + k\nabla_YZ, \quad
\nabla_X(hY + kZ) = h\nabla_XY + (Xh)Y + k\nabla_XZ + (Xk)Z, 
$$
for any vector fields $X, Y, Z$ and functions $h, k$ on $M$. 

For a system of local coordinates $x^\lambda\, (\lambda = 1, 2, \dots, m)$, we write 
$$
\nabla_{\frac{\pa}{\pa x^\mu}}{\textstyle\frac{\pa}{\pa x^\nu}} \ = \ \Gamma^\lambda_{\mu \nu}\ {\textstyle\frac{\pa}{\pa x^\lambda}}, 
$$
using the Christoffel symbols (coefficients of the connection) 
$\Gamma^\lambda_{\mu \nu}$ and the Einstein convention. 

In general, for a given mapping $g : N \to (M, \nabla)$, we define the notion of 
{\it covariant derivative} $\nabla_\eta^g v : N \to TM$ 
of a vector field $v : N \to TM$ along $g$ 
by a vector field $\eta : N \to TN$ over a manifold $N$ (see \cite{FSUY}). 
Using a local presentation 
$\eta(t) = \eta^i(t)\frac{\pa}{\pa t^i}(t)$, 
$v(t) = v^\lambda(t)\frac{\pa}{\pa x^\lambda}(g(t))$,  
for local coordinates $t = (t_1, \dots, t_n)$ of $N$ and $x = (x_1, \dots, x_m)$ of $M$, 
and using the Einstein convention, we define 
$$
(\nabla_\eta^g v)(t) := \left\{ \eta^i(t) \frac{\pa v^\lambda }{\pa t^i}(t)
+ \Gamma^\lambda_{\mu\nu}(g(t))\eta^i(t) \frac{\pa g^\mu}{\pa t^i}(t)v^\nu(t)\right\}\frac{\pa}{\pa x^\lambda}(g(t)). 
$$
Here $g^\mu = x^\mu\circ g$. The definition is naturally derived from the parallelisms on $M$ induced by the connection 
$\nabla$. 

Then we have, as in the usual case, 
$$
\nabla_{h\eta + k\xi}^f v = h\nabla_\eta^f v + k \nabla_\xi^f v, 
\quad
\nabla_\eta^f(hv + kw) = h\nabla_\eta^f v + k\nabla_\eta^f w 
+ (\eta h)v + (\eta k)w, 
$$
for any vector fields $v, w$ along $f$, $\eta, \xi$ over $N$, and functions $h, k$ on $N$. 

If $g : M \to M$ is the identity mapping, then $\nabla_\eta^g v$ is just 
the ordinary covariant derivative $\nabla_\eta v$ for vector fields $\eta, v$ over $M$. 
The covariant derivative along a mapping is well-defined also for any tensor field over the mapping, 
that is compatible with any contractions. 

Then we get the notion of geodesics: 
A curve $\varphi : I \to M, \varphi = \varphi(s)$ is called a $\nabla$-{\it geodesic}  
if $\nabla_{\pa/\pa s}^\varphi \frac{d \varphi}{ds} = 0$. 
For $x \in M, v \in T_xM$, let $\varphi(x, v, s)$ denote the $\nabla$-geodesic 
determined by the differential equation
$$
\dfrac{\pa^2\varphi}{\pa s^2}^\lambda(x, v, s) \ + \ \Gamma^\lambda_{\mu \nu}(\varphi(x, v, s))
\dfrac{\pa \varphi}{\pa s}^\mu(x, v, s)\dfrac{\pa \varphi}{\pa s}^\nu(x, v, s) = 0, 
$$
with the initial conditions $\varphi(x, v, 0) = x$ and $\frac{\pa \varphi}{\pa s}(x, v, 0) = v$. 
Here $\varphi^\lambda$ denotes the $\lambda$-th component $x^\lambda\circ \varphi$ of $\varphi$. 
Note that $\varphi : U (\subset TM\times\R) \to M$ 
is defined on an open neighborhood $U$ of $TM\times\{ 0\}$. 

\bel
\label{geodesic-expression}
At each point of M, there exists an open neighborhood $U$ of the point such that 
the $\nabla$-geodesics $\varphi(x, v, s)$ are written 
$$
\varphi(x, v, s) = x + s\, v + \frac{1}{2}s^2\, h(x, v, s)
$$
for some $C^\infty$ mapping $h(x, v, s)$ on an open neighborhood of $TU\times \{ 0\}$ in 
$TU\times\R$. 
\enl

\Proof
Since $\varphi(x, v, 0) = x$, we have locally 
$\varphi(x, v, s) = x + s\tilde{\varphi}(x, v, s)$ for some $C^\infty$ mapping $\tilde{\varphi}$. 
Then $\frac{\pa \varphi}{\pa s}(x, v, s) = \tilde{\varphi}(x, v, s) + s\frac{\pa \tilde{\varphi}}{\pa s}(x, v, s)$. 
Since $\tilde{\varphi}(x, v, 0) = \frac{\pa \varphi}{\pa s}(x, v, 0) = v$, we have locally 
$\tilde{\varphi}(x, v, s) = v + \frac{1}{2}s\, h(x, v, s)$ for some $C^\infty$ mapping $h$. 
Thus we have the result. 
\QED

\bel
\label{reminder}
Let $\varphi(x, v, s) = x + s\, v + \frac{1}{2}s^2\, h(x, v, s)$ be the expression of the $\nabla$-geodesics 
as in Lemma \ref{geodesic-expression}. Then we have
\begin{align*}
h^\lambda(x, v, 0) & = \frac{\pa^2 \varphi^\lambda}{\pa s^2}(x, v, 0) = - \Gamma^\lambda_{\mu\nu}(x)v^\mu v^\nu, 
\\
\frac{\pa h^\lambda}{\pa x^\kappa}(x, v, 0) & = - \Gamma^\lambda_{\mu\nu,\kappa}(x)v^\mu v^\nu, 
\\
\frac{\pa h^\lambda}{\pa v^\rho}(x, v, 0) & = - \Gamma^\lambda_{\rho\nu}(x)v^\nu - \Gamma^\lambda_{\mu\rho}(x)v^\mu, 
\\
\frac{\pa h^\lambda}{\pa s}(x, v, 0) & = \frac{1}{3}\frac{\pa^3\varphi}{\pa s^3}(x, v, 0) 
= \frac{1}{3}(-\Gamma_{\mu\nu,\kappa}^\lambda + \Gamma_{\rho\kappa}^\lambda\Gamma_{\mu\nu}^\rho 
+ \Gamma_{\kappa\rho}^\lambda\Gamma_{\mu\nu}^\rho)v^\mu v^\nu v^\kappa
\end{align*}
where $\Gamma^\lambda_{\mu\nu,\kappa} = \pa \Gamma^\lambda_{\mu\nu}/ \pa x^\kappa$. 
\enl

\Proof
For a fixed $(x, v)$,  $\varphi(x, v, s)$ is a $\nabla$-geodesic, therefore we have 
$$
\dfrac{\pa^2 \varphi^\lambda}{\pa s^2}(x, v, s) 
= - \Gamma^\lambda_{\mu\nu}(\varphi(x, v, s))\frac{\pa \varphi}{\pa s}^\mu(x, v, s)\frac{\pa \varphi}{\pa s}^\nu(x, v, s). 
$$
Since $\frac{\pa \varphi}{\pa s}(x, v, 0) = v$, we have that $\frac{\pa^2 \varphi}{\pa s^2}(x, v, 0) = h(x, v, 0)$, 
and that 
$h^\lambda(x, v, 0) = - \Gamma^\lambda_{\mu\nu}(x)v^\mu v^\nu$. 
Moreover we have 
\begin{align*}
\dfrac{\pa^3 \varphi^\lambda}{\pa s^3}(x, v, s) 
= & - \Gamma^\lambda_{\mu\nu,\kappa}(\varphi(x, v, s))\frac{\pa \varphi}{\pa s}^\kappa(x, v, s)\frac{\pa \varphi}{\pa s}^\mu(x, v, s)\frac{\pa \varphi}{\pa s}^\nu(x, v, s) 
\\
& 
- \Gamma^\lambda_{\mu\nu}\frac{\pa^2 \varphi}{\pa s^2}^\mu(x, v, s)\frac{\pa \varphi}{\pa s}^\nu(x, v, s)
- \Gamma^\lambda_{\mu\nu}\frac{\pa \varphi}{\pa s}^\mu(x, v, s)\frac{\pa^2 \varphi}{\pa s^2}^\nu(x, v, s). 
\end{align*}
By setting $s = 0$ and by the first equality, we have the fourth equality. 
By differentiating both sides of 
$h^\lambda(x, v, 0) = - \Gamma^\lambda_{\mu\nu}(x)v^\mu v^\nu$ (the first equality), 
by $x^\kappa$ and $v^\rho$ we have the second and the third equalities. 
\QED

\

A connection $\nabla$ on $M$ is called {\it torsion free} if 
$T(X, Y) := \nabla_XY - \nabla_YX - [X, Y] = 0, $ for any vector fields $X, Y$ on $M$. A connection 
$\nabla$ is torsion free if and only if, for any system of local coordinates, 
$\Gamma^\lambda_{\mu\nu} = \Gamma^\lambda_{\nu\mu}, 1 \leq \lambda, \mu, \nu \leq m$. 

\ber
\label{torsion-free}
{\rm 
We observe that the equation on geodesics 
$$
\dfrac{\pa^2\varphi}{\pa s^2}^\lambda(x, v, s) \ + \ \Gamma^\lambda_{\mu \nu}(\varphi(x, v, s))
\dfrac{\pa \varphi}{\pa s}^\mu(x, v, s)\dfrac{\pa \varphi}{\pa s}^\nu(x, v, s) = 0, 
$$
is symmetric on the indices $\mu, \nu$. Therefore 
the geodesics $\varphi(x, v, s)$ and the tangent surfaces $\Tan(\gamma)$ 
remain same if the connection $\Gamma^\lambda_{\mu \nu}$ is replaced by the torsion free 
connection $\frac{1}{2}(\Gamma^\lambda_{\mu \nu} + \Gamma^\lambda_{\nu \mu})$, 
in other word, if $\nabla$ is replaced by the torsion free connection $\widetilde{\nabla}$, 
defined by $\widetilde{\nabla}_XY = \nabla_XY - \frac{1}{2}T(X, Y)$. 
}
\enr

\

Let $\gamma : I \to M$ be a curve which is not necessarily a geodesic nor 
an immersed curve. Then the first derivative $(\nabla\gamma)(t)$ means just the velocity 
vector field $\gamma'(t)$. The second derivative $(\nabla^2\gamma)(t)$ is defined, in terms of 
covariant derivative along the curve $\gamma$, by 
$$
(\nabla^2\gamma)(t) := \nabla_{\pa/\pa t}^\gamma (\nabla\gamma)(t). 
$$
Note that $\gamma$ is a $\nabla$-geodesic if and only if $\nabla^2 \gamma = 0$. 
In general, we define $k$-th covariant derivative of $\gamma$ inductively by 
$$
(\nabla^k\gamma)(t) := \nabla_{\pa/\pa t}^\gamma (\nabla^{k-1}\gamma)(t), \ (k \geq 2). 
$$
Then we have by direct calculations: 

\bel
\label{iteration}
\begin{align*}
(\nabla\gamma)^\lambda & = (\gamma')^\lambda, 
\\
(\nabla^2\gamma)^\lambda & = (\gamma'')^\lambda + \Gamma^\lambda_{\mu\nu}(\gamma')^\mu(\gamma')^\nu, 
\\
(\nabla^3\gamma)^\lambda & 
= (\gamma''')^\lambda + 
(\Gamma^\lambda_{\mu\nu, \kappa} + \Gamma^\lambda_{\kappa\rho}\Gamma^\rho_{\mu\nu})
(\gamma')^\mu(\gamma')^\nu(\gamma')^\kappa 
+ (2\Gamma^\lambda_{\mu\nu} + \Gamma^\lambda_{\nu\mu})(\gamma')^\mu(\gamma'')^\nu. 
\end{align*}
\enl

\section{Tangent surface and frontal}
\label{Tangent surface and frontal}

Let $\gamma : I \to M$ be a $C^\infty$ immersion from an open interval $I$. 
Then the tangent surface to $\gamma$ is the ruled surface by tangent $\nabla$-geodesics to $\gamma$. 
More precisely the $\nabla$-tangent surface $f = \nabla{\mbox{\rm -}}\Tan(\gamma) : V (\subset I \times \R) \to M$ to the curve 
$\gamma$ is defined by 
$$
f (t, s) := \varphi(\gamma(t), \gamma'(t), s), 
$$
on an open neighborhood $V$ of $I \times \{0\}$. 

The mapping $f$ has singularity at least along $\{ s = 0\}$. In fact 
$f_*(\frac{\pa}{\pa t})(t, 0) = \gamma'(t) = f_*(\frac{\pa}{\pa s})(t, 0)$ and the kernel of the differential $f_*$ 
is generated by the vector field $\eta = \frac{\pa}{\pa t} - \frac{\pa}{\pa s}$ along $\{ s = 0 \}$ on the $t$-$s$-plane. 

A map-germ $f : (\R^n, p) \to M$, $n \leq m = \dim(M)$ is called a {\it frontal} if 
there exists a $C^\infty$ {\it integral} lifting $\widetilde{f} : (\R^n, p) \to \Gr(n, TM)$ of $f$. 
Here $\Gr(n, TM)$ means the Grassmannian bundle over $M$ consisting of 
tangential $n$-planes in $TM$, 
$$
\Gr(n, TM) := \{ \Pi \mid \Pi \subseteq T_xM, \Pi {\mbox{\rm \ is a linear subspace,\ }}\dim(\Pi) = n, \, x \in M\}, 
$$ 
with the canonical projection $\pi : \Gr(n, TM) \to M$ to the base manifold $M$, 
and we call $\tilde{f}$ is {\it integral} if $f_*(T_q\R^n) \subseteq \widetilde{f}(q)$ for any $q$ 
in a neighborhood of $p$ in 
$\R^n$, after taking a representative of $f$. 
The definition generalizes the preceding definition of \lq\lq frontal" in the case $m = n+1$ (see \cite{FSUY}). 
The definitions in the case $m = n+1$ are equivalent to each other as is easily seen. 
Note that, in \cite{Ishikawa4}, 
we have introduced the same notion of frontal mapping under the restriction that the locus of immersive points 
of $f$ is dense, where the integral lifting $\widetilde{f}$ is uniquely determined. 

Let $f : (\R^n, p) \to M$ is a frontal and $\widetilde{f}$ is an integral lifting of $f$. 
Then there exists a frame $V_1, V_2, \dots, V_n : (\R^n, p) \to TM$ 
along $f$ associated with $\widetilde{f}$ such that 
$$
\widetilde{f}(q) = \langle V_1(q), V_2(q), \dots, V_n(q)\rangle_{\R}, 
$$
for any $q$ in a neighborhood of $p$ in 
$\R^n$. 
Then there is a $C^\infty$ function-germ $\sigma : (\R^n, p) \to \R$ such that 
$$
(\dfrac{\pa f}{\pa t_1} \wedge \dfrac{\pa f}{\pa t_2} \wedge \cdots \wedge \dfrac{\pa f}{\pa t_n})(t) 
= \sigma(t)(V_1 \wedge V_2 \wedge \cdots \wedge V_n)(t), 
$$
as germs of $n$-vector fields $(\R^n, p) \to \wedge^n TM$ over $f$. 
Then the singular locus (non-immersive locus) $S(f)$ of $f$ coincides with the zero locus $\{ \sigma = 0\}$ 
of $\sigma$. 
We call $\sigma$ a {\it signed area density function} or briefly an {\it $s$-function} of the frontal $f$. 
Note that the function $\sigma$ is essentially the same thing with the function $\lambda$ introduced in 
\cite{KRSUY}\cite{FSUY} in the case $\dim(M) = 3$. However we avoid the notation $\lambda$ here because 
we use it for index of Christoffel symbols. 

We say that a frontal $f : (\R^n, p) \to M$ has a {\it non-degenerate} singular point at $p$ if 
the $s$-function $\sigma$ of $f$ satisfies $\sigma(p) = 0$ and $d\sigma(p) \not= 0$. 
The condition is independent of the choice of the integral lifting $\widetilde{f}$ 
and the associated frame $V_1, V_2, \dots, V_n$. 
If $f$ has a non-degenerate singular point at $p$, then $f$ is of corank $1$ such that the singular locus 
$S(f) \subset (\R^n, p)$ is a regular hypersurface. 

In this paper we concern with only the cases $n = 1$ and $n = 2$. 

\

Returning to our situation, we have 

\bel
\label{non-degenerate-frontal}
Suppose $(\nabla\gamma)(t_0)$ and $(\nabla^2\gamma)(t_0)$ are linearly independent. 
Then the germ of tangent surface $\nabla{\mbox{\rm -}}\Tan(\gamma)$ 
is a frontal with the non-degenerate singular point at $(t_0, 0)$ and 
with the singular locus $S(\nabla{\mbox{\rm -}}\Tan(\gamma)) = \{ s = 0\}$. 
\enl

\Proof
We set $f(t, s) = \Tan(\gamma)(t, s) = \varphi(\gamma(t), \gamma'(t), s)$. 
By Lemma \ref{geodesic-expression}, write $\varphi = x + s\, v + \frac{1}{2}s^2\, h$. Then we have 
$$
f(t, s) = \gamma(t) + s\, \gamma'(t) + \frac{1}{2} s^2 \, h(\gamma(t), \gamma'(t), s), 
$$
and 
\begin{align*}
\frac{\pa f}{\pa t} & =  \gamma' + s\gamma'' + 
\frac{1}{2}s^2\, (\gamma')^\mu\frac{\pa h}{\pa x^\mu}(\gamma, \gamma', s) + 
\frac{1}{2}s^2\, (\gamma'')^\nu\frac{\pa h}{\pa v^\nu}(\gamma, \gamma', s), 
\\
\frac{\pa f}{\pa s} & =  \gamma' + s\, h(\gamma, \gamma', s) + \frac{1}{2}s^2\, \frac{\pa h}{\pa s}(\gamma, \gamma', s). 
\end{align*}
Then we see that $S(f) \supseteq \{ s = 0 \}$. 

Let $s \not= 0$. Then
\begin{equation*}
\begin{split}
\frac{1}{s}(\frac{\pa f}{\pa t} - \frac{\pa f}{\pa s}) = 
& \ 
\gamma'' + \frac{1}{2}s\, (\gamma')^\mu\frac{\pa h}{\pa x^\mu}(\gamma, \gamma', s) 
+ \frac{1}{2}s\, (\gamma'')^\nu\frac{\pa h}{\pa v^\nu}(\gamma, \gamma', s) 
\\
& 
\hspace{0.5truecm} 
 - h(\gamma, \gamma', s) - \frac{1}{2}s\, \frac{\pa h}{\pa s}(\gamma, \gamma', s). 
\end{split}
\end{equation*}

We define $F(t, s)$ by the right hand side. Then  $F(t, s) = \frac{1}{s}(\frac{\pa f}{\pa t} - \frac{\pa f}{\pa s})$ if 
$s \not= 0$. Moreover $F$ is $C^\infty$ also on $s = 0$ and 
$$
F(t, 0) = \gamma''(t) - h(\gamma(t), \gamma'(t), 0). 
$$
By Lemma \ref{reminder}, we have 
$$
h^\lambda(\gamma(t), \gamma'(t), 0) = - \Gamma^\lambda_{\mu\nu}(\gamma(t))\, (\gamma'(t))^\mu(\gamma'(t))^\nu. 
$$
Hence we have 
$$
F^\lambda(t, 0) = (\gamma''(t))^\lambda + \Gamma^\lambda_{\mu\nu}(\gamma(t))\, (\gamma'(t))^\mu(\gamma'(t))^\nu 
= (\nabla^2\gamma)^\lambda(t). 
$$
Therefore if 
$(\nabla \gamma)(t), (\nabla^2\gamma)(t)$ are linearly independent at $t = t_0$, then 
$\frac{\pa f}{\pa t}(t, s)$ and $F(t, s)$ are linearly independent around $(t_0, 0)$ and satisfies that 
$$
(\frac{\pa f}{\pa t} \wedge \frac{\pa f}{\pa s})(t, s) = -s (\frac{\pa f}{\pa t} \wedge F)(t, s). 
$$
Therefore we see that $\frac{\pa f}{\pa t}(t, s)$ and $F(t, s)$ define the integral lifting of $f$, 
$f$ is frontal with non-degenerate singular point at $(t_0, 0)$, and that 
$S(f) = \{ s = 0\}$.
\QED

\ber
{\rm 
In the Euclidean case, 
Lemma \ref{non-degenerate-frontal} holds globally on $I\times \R$. However, even in a locally projectively flat case, 
Lemma \ref{non-degenerate-frontal} holds just locally near $I \times \{ 0\}$. 
For example, let $M$ be the standard three dimensional sphere $S^3 \subset \R^4$ 
with the standard (Levi-Civita) connection. Then geodesics in $S^3$ are given by great circles (with periodic parametrizations) and 
we observe, via the natural double covering $S^3 \to \R P^3$, 
that the tangent surface to any curve in $S^3$ 
has singularities not only along the original curve, but also along the antipodal of the curve (cf. \cite{INS}). 
}
\enr

\section{Cuspidal edge and folded umbrella}
\label{Cuspidal edge and folded umbrella}

Let $f : (\R^2, p) \to M^3$ be a frontal with a non-degenerate singular point at $p$ 
and $\widetilde{f} : (\R^2, p) \to \Gr(2, TM)$ the integral lifting of $f$. 
Let $V_1, V_2 : (\R^2, p) \to TM$ be an associated frame with $\widetilde{f}$. 
Let $L : (\R^2, p) \to T^*M \setminus \zeta$ be an annihilator of 
$\widetilde{f}$. The condition is that $\langle L, V_1\rangle = 0, \langle L, V_2\rangle = 0$. 
Here $\zeta$ means the zero section. 
Let $c : (\R, t_0) \to (\R^2, p)$ be a parametrization of the singular locus $S(f)$, $p = c(t_0)$,  and 
$\eta : (\R^2, p) \to T\R^2$ be a vector field which restricts to the kernel field of $f_*$ on $S(f)$. 
Suppose that $V_2(p) \not\in f_*(T_p\R^2)$. 
Then, for any affine connection $\nabla$ on $M$, we define 
$$
\psi(t) := \langle L(c(t)),  (\nabla_\eta^f V_2)(c(t))\rangle. 
$$
Note that the vector field $(\nabla_\eta^f V_2)(c(t))$ is independent of the extension $\eta$ and the choice of affine connection $\nabla$, since $\eta\vert_{S(f)}$ is a kernel field of $f_*$. 
We call the function $\psi(t)$
the {\it characteristic function} of $f$. 

Then the following characterizations of cuspidal edges and folded umbrellas are given in \cite{KRSUY}\cite{FSUY}: 

\bet
\label{characterization}
{\rm (Theorem 1.4 of \cite{FSUY})}. 
Let $f : (\R^2, p) \to M^3$ be a germ of frontal with a non-degenerate singular point at $p$. 
Let $c : (\R, t_0) \to (\R^2, p)$ be a parametrization of the singular locus of $f$. 
Suppose $f_*c'(t_0) \not= 0$. Then, for the characteristic function $\psi$, 
\\
{\rm (1)} 
$f$ is diffeomorphic to the cuspidal edge if and only if $\psi(t_0) \not= 0$. 
\\
{\rm (2)} $f$ is diffeomorphic to the folded umbrella if and only if
$\psi(t_0) = 0, \psi'(t_0) \not= 0$. 
\ent

Note that the conditions appeared in Theorem \ref{characterization} are invariant under diffeomorphism equivalence.

\ber
{\rm 
In the situation of Theorem \ref{characterization}, we set $\gamma(t) = f(c(t))$. Then we have 
$\psi(t_0) \not= 0$ if and only if $V_1(c(t_0)), V_2(c(t_0)), (\nabla_\eta^f V_2)(c(t_0))$ are linearly independent. 
}
\enr

\

The above construction is generalized to the case $m = \dim(M) \geq 4$. 
In general, let $f : (\R^2, p) \to M^m, m \geq 4,$ be a frontal with a non-degenerate singular point at $p$ 
and $\widetilde{f} : (\R^2, p) \to \Gr(2, TM)$ the integral lifting of $f$. 
Let $V_1, V_2 : (\R^2, p) \to TM$ be an associated frame with $\widetilde{f}$. 
We take a coframe $L_1, \dots, L_{m-2} : (\R^2, p) \to T^*M$ satisfying that 
$$
\langle L_i(t, s), V_1(t, s)\rangle = 0, \langle L_i(t, s), V_2(t, s)\rangle = 0, \quad (1 \leq i \leq m-2), 
$$
and that $L_1(t, s), \dots, L_{m-2}(t, s)$ are linearly independent for $(t, s) \in (\R^2, p)$. 
We define the {\it characteristic (vector valued) function} $\psi : (\R, t_0) \to \R^{m-2}$ by $\psi = (\psi_1, \dots, \psi_{m-2})$, 
$$
\psi_i(t) := \langle L_i(c(t)),  (\nabla_\eta^f V_2)(c(t))\rangle. 
$$

Let $g : (\R^n, p) \to (\R^\ell, q)$ be a map-germ. 
A map germ $f : (\R^n, p) \to \R^{\ell+r}$ is called an {\it opening} of $g$ if there exist
functions $h_1, \dots, h_r : (\R^n, p) \to \R$ and functions 
$a_{ij} : (\R^n, p) \to \R, (1 \leq i \leq r, 1 \leq j \leq \ell)$ such that 
$$
dh_i = \sum_{j=1}^\ell a_{ij}dg_j,  \quad f = (g, h_1, \dots, h_r), 
$$
(see for example \cite{Ishikawa4}). If $\ell = n$, then the condition on $h$ is equivalent to that 
$f$ is frontal associated with an integral lifting $\widetilde{f} : (\R^n, p) 
\to \Gr(n, T\R^{n+r})$ having Grassmannian coordinates $(a_{ij})$  
such that $\widetilde{f}(p)$ projects isomorphically to 
$T_{g(p)}\R^n$ by the projection $\R^{n+r} = \R^n\times \R^r \to \R^n$. 

\ber
{\rm 
A map-germ $f : (\R^n, p) \to M$ is a frontal with a non-degenerate singular point at $p$ 
if and only if $f$ is diffeomorphic to an opening of a map-germ $g : (\R^n, p) \to (\R^n, q)$ 
of Thom-Boardman singularity type $\Sigma^1$ at $p$, i.e. 
$g$ is of corank one and $j^1g : (\R^n, p) \to J^1(\R^n, \R^n)$ is transversal to 
the variety of singular $1$-jets (see for example \cite{GG}). 
}
\enr

We can summarize several known results as those on openings of the fold: 

\bet
\label{characterization2}
Let $f : (\R^2, p) \to M^m, m \geq 2$ be a germ of frontal with a non-degenerate singular point at $p$, 
$\widetilde{f} : (\R^2, p) \to \Gr(2, TM)$ the integral lifting of $f$ and 
$V_1, V_2 : (\R^2, p) \to TM$ an associated frame with $\widetilde{f}$. 
Let $c : (\R, t_0) \to (\R^2, p)$ be a parametrization of the singular locus of $f$. 
Suppose $f_*c'(t_0) \not= 0$. 
Then $f$ is diffeomorphic to an opening of the fold, namely to the germ $(u, w) \mapsto (u, \frac{1}{2}w^2)$. 
Moreover we have: 

{\rm (0)} Let $m = 2$. Then $f$ is diffeomorphic to the fold. 

{\rm (1)} Let $m \geq 3$. Then $f$ is diffeomorphic to the cuspidal edge if and only if
$\psi(t_0) \not= 0$. 

{\rm (2)} Let $m = 3$. Then $f$ is diffeomorphic to the folded umbrella if and only if
$\psi(t_0) = 0, \psi'(t_0) \not= 0$. 
\ent

\

\noindent
{\it Proof of Theorem \ref{characterization2}:} 
The assertion (0) follows from Whitney's theorem (see \cite{Whitney}\cite{SUY}\cite{Saji}). 
(1) The condition $\psi(t_0) \not= 0$ is equivalent to that $f$ is a front, namely, that $\widetilde{f}$ is an immersion. 
Suppose $\dim(M) = 3$. Then by Proposition 1.3 of \cite{KRSUY}, 
we see that $f$ is diffeomorphic to the cuspidal edge. 
In general cases $m \geq 2$, we see that 
there exists a submersion $\pi : (M, f(p)) \to (\R^2, 0)$ such that 
$$
T_{f(p)}(\pi^{-1}(0)) + \widetilde{f}(p) = T_{f(p)}M, 
$$
and that $\pi\circ f$ satisfies the same condition with $f$, 
i.e., $\pi\circ f$ is a frontal with the non-degenerate singular point at $p$ with the same singular locus with 
$f$ and $\eta(c(t_0))$ and $c'(t_0)$ are linearly independent, but $m = 2$. Thus by the assertion (0), 
the map-germ $\pi\circ f$ is diffeomorphic to a fold. 
Moreover we see $f$ is an opening of $\pi\circ f$ because 
$f$ is frontal. 
The condition that $\widetilde{f}$ is an immersion is equivalent, in this case, 
to that $f$ is a versal opening of $\pi\circ f$ (\S 6 of \cite{Ishikawa4}). 
In fact, up to diffeomorphism equivalence, let 
$$
f(u, w) = (u, \frac{1}{2}w^2, h_1(u, w), \dots, h_r(u, w)), 
$$
$(m = 2+r)$ and
$$
dh_i(u, w) = a_i(u, w)du + b_i(u, w)d(\frac{1}{2}w^2) = a_i(u, w)du + wb_i(u, w)dw, 
$$
for some functions $a_i, b_i, 1 \leq i \leq r$. 
Then we have 
$$
\frac{\pa f}{\pa u} = (1, 0, a_1, \dots, a_r), \quad \frac{\pa f}{\pa w} = (0, w, w b_1, \dots, w b_r)
$$
and the pair $V_1 = \frac{\pa f}{\pa u}, V_2 = \frac{1}{w}\frac{\pa f}{\pa w} = (0, 1, b_1, \dots, b_r)$
gives a frame of the frontal $f$. Moreover $\eta = \frac{\pa}{\pa w}$ gives the kernel field of $f_*$ 
along $\{w = 0\}$. 
For any connection $\nabla$, we have the characteristic vector field 
$\nabla^f_\eta V_2(u, 0) = (0, 0, \frac{\pa b_1}{\pa w}, \dots, \frac{\pa b_r}{\pa w})$. 
Then $\psi(0) \not= 0$ if and only if $V_1(0, 0), V_2(0, 0), (\nabla^f_\eta V_2)(0, 0)$ 
are linearly independent. The condition is equivalent to that $h_1(0, w), \dots, h_r(0, w)$ 
generate $m_1^3/m_1^4$ over $\R$ and it is equivalent to that $f$ is a versal opening of 
the Whitney's cusp. Here $m_1$ means the ideal consisting of function-germs $h(w)$ with $h(0) = 0$. 
Then we see that $f$ is diffeomorphic to cuspidal edge (see Proposition 6.8 (3) $\ell = 2$ of \cite{Ishikawa4}). 
The assertion (2) follows from Theorem 1.4 of \cite{FSUY}. 
\QED

\section{Characteristic vector field}
\label{Characteristic vector field}

Let $\gamma : I \to M$ be an immersion. We set $f = \nabla{\mbox{\rm -}}\Tan(\gamma)$ and 
suppose $\nabla\gamma, \nabla^2\gamma$ are linearly independent at $t = t_0$. 
Note that $\eta = \frac{\pa}{\pa t} - \frac{\pa}{\pa s}$ generates the field of kernels $\Ker(f_*)$ of the differential 
$f_*$ along $s = 0$. 
Let $F(t, s) = \frac{1}{s}(\frac{\pa f}{\pa t} - \frac{\pa f}{\pa s})$ as in the proof of Lemma \ref{non-degenerate-frontal}. 
Then we have
$$
(\nabla_\eta^f F)(t, s) = \{(\frac{\pa}{\pa t} - \frac{\pa}{\pa s})F^\lambda + (\Gamma^\lambda_{\mu\nu}(f(t, s)) 
((\frac{\pa}{\pa t} - \frac{\pa}{\pa s}) f^\mu) F^\nu \}\frac{\pa}{\pa x^\lambda}(f(t, s)). 
$$
Therefore we have 
$$
(\nabla_\eta^f F)(t, 0) = (\frac{\pa F}{\pa t} - \frac{\pa F}{\pa s})(t, 0). 
$$

We call the vector field $(\nabla_\eta^f F)(t, 0)$ along $\gamma$ the {\it characteristic vector field} of 
$\nabla{\mbox{\rm -}}\Tan(\gamma)$. 

By straightforward calculations we have 

\bel
\label{characteristic vector}
The characteristic vector field of $\nabla{\mbox{\rm -}}\Tan(\gamma)$ is given by 
\begin{align*}
(\nabla_\eta^f F)^\lambda(t, 0) 
& = 
(\gamma''')^\lambda + (\Gamma^\lambda_{\mu\nu,\kappa} + \frac{1}{2}\Gamma^\lambda_{\rho\mu}\Gamma^\rho_{\nu\kappa}
+ \frac{1}{2}\Gamma^\lambda_{\mu\rho}\Gamma^\rho_{\nu\kappa})(\gamma')^\mu(\gamma')^\nu(\gamma')^\kappa
+ \frac{3}{2}(\Gamma^\lambda_{\mu\nu} + \Gamma^\lambda_{\nu\mu})(\gamma')^\mu(\gamma'')^\nu. 
\end{align*}
\enl

\bel
\label{third derivative}
$(\nabla_\eta^f F)(t, 0) = (\nabla^3 \gamma)(t)$ if the affine connection $\nabla$ is torsion free. 
\enl

\Proof We compare Lemma \ref{characteristic vector} and Lemma \ref{iteration}. 
The equality $\frac{3}{2}(\Gamma^\lambda_{\mu\nu} + \Gamma^\lambda_{\nu\mu}) = 
2\Gamma^\lambda_{\mu\nu} + \Gamma^\lambda_{\nu\mu}$ holds if and only if 
$\Gamma^\lambda_{\mu\nu} = \Gamma^\lambda_{\nu\mu}$. Then 
the equality 
$\frac{1}{2}\Gamma^\lambda_{\rho\mu}\Gamma^\rho_{\nu\kappa}
+ \frac{1}{2}\Gamma^\lambda_{\mu\rho}\Gamma^\rho_{\nu\kappa} 
= \Gamma^\lambda_{\mu\rho}\Gamma^\rho_{\nu\kappa}$ holds. 
\QED

\section{Euclidean case}
\label{Euclidean case}

Now we show Theorem \ref{characterization-theorem} in Euclidean case, using the characterization results 
of singularities prepared in previous sections. 

Let $M = {E}^m$ be the Euclidean $m$-space with the connection $\Gamma^\lambda_{\mu\nu} = 0$. 
Then $\varphi(x, v, s) = x + s\, v$, that is $h(x, v, s) = 0$. 
If $\gamma'(t_0) \not= 0$, then the tangent surface is given by
$f(t, s) = \nabla{\mbox{\rm -}}\Tan(\gamma)(t, s) = \gamma(t) + s\gamma'(t)$ near $t_0 \times \R$. 
Then $\frac{\pa f}{\pa t} = \gamma'(t) + s\gamma''(t), 
\frac{\pa f}{\pa s} = \gamma'(t)$ and we have $\frac{\pa f}{\pa t} - \frac{\pa f}{\pa s} = s\gamma''(t)$. 
$F(t, s) = \frac{1}{s}(\frac{\pa f}{\pa t} - \frac{\pa f}{\pa s}) = \gamma''(t)$. Moreover we have 
$(\nabla_{\pa/\pa t - \pa/\pa s}^f F)(t, 0) = \gamma'''(t)$. 
Suppose $\gamma'(t_0)$ and $\gamma''(t_0)$ are linearly independent. 
By Lemma \ref{non-degenerate-frontal}, $f$ is a frontal with non-degenerate singular point at $(t_0, 0)$. 
We apply Theorem \ref{characterization} to this situation. We take $\eta = \frac{\pa}{\pa t} - \frac{\pa}{\pa s}$ and 
$c(t) = (t, 0)$. Then $c'(t_0)$ and $\eta(c(t_0))$ are linearly independent, that is, $f_*(c'(t_0)) \not= 0$. 
Take a coframe, namely, a system of 
germs of $1$-forms $L_i = L_i(t, s) : (I\times \R, (t_0, 0)) \to T^*{E}^m \setminus \zeta$ along $f$, $i = 1, \dots, m-2$, which satisfy that 
$$
\langle L_i(t, s), \gamma'(t)\rangle = 0, \quad \langle L_i(t, s), \gamma''(t)\rangle = 0, (i = 1, \dots, m-2), 
$$
and $L_1(t, s), \dots L_{m-2}(t, s)$ are linearly independent for $(t, s) \in (I\times\R, (t_0, 0))$. 
Here $\zeta$ means the zero section. Actually we can take $L$ to be independent of $s$ in this case. 
Set $\ell_i(t) = L_i(t, 0)$ and 
set $\psi_i(t) = \langle \ell_i(t), \gamma'''(t)\rangle$ and 
$$
\psi(t) = (\psi_1(t), \dots, \psi_{n-2}(t)) = (\langle \ell_1(t), \gamma'''(t)\rangle, \dots, \langle \ell_{m-2}(t), \gamma'''(t)\rangle), 
$$
the characteristic function.  Then we have that 
$\psi(t_0) = 0$ if and only if $\gamma'(t_0), \gamma''(t_0), \gamma'''(t_0)$ are linearly dependent. 
We have 
$$
\langle \ell_i(t), \gamma'(t)\rangle = 0, \quad \langle \ell_i(t), \gamma''(t)\rangle = 0, (i = 1, \dots, m-2), 
$$
and 
\begin{align*}
0 & = \langle \ell_i'(t), \gamma'(t)\rangle + \langle \ell_i(t), \gamma''(t)\rangle = \langle \ell_i'(t), \gamma'(t)\rangle, 
\\
0 & = \langle \ell_i'(t), \gamma''(t)\rangle + \psi_i(t). 
\end{align*}
Suppose $\psi(t_0) = 0$. Then $\langle \ell_i'(t_0), \gamma''(t_0)\rangle = 0$ for any $i, (1 \leq i \leq m-2)$. 
Since we have also $\langle \ell_i'(t_0), \gamma'(t_0)\rangle = 0$, we obtain $\langle \ell_i'(t_0), \gamma'''(t_0)\rangle = 0$, 
because $\gamma'(t_0)$ and $\gamma''(t_0)$ are linearly independent and 
$\gamma'(t_0), \gamma''(t_0), \gamma'''(t_0)$ are linearly dependent. 
Now we have 
$$
\psi_i'(t) = \langle \ell_i'(t), \gamma'''(t)\rangle + \langle \ell_i(t), \gamma^{(4)}(t)\rangle.
$$
Therefore we have $\psi_i'(t_0) = \langle \ell_i(t_0), \gamma^{(4)}(t_0)\rangle$. 
Thus under the condition $\psi(t_0) = 0$, we have 
$\psi'(t_0) \not= 0$ if and only if $\gamma'(t_0), \gamma''(t_0), \gamma^{(4)}(t_0)$ are linearly independent. 

Thus Theorem \ref{characterization} and Theorem \ref{characterization2} (1) imply 
Theorem \ref{characterization-theorem} (1)(2) in the Euclidean case.

\section{Proof of the characterization theorem}
\label{Proof of the characterization theorem}

\noindent
{\it Proof of Theorem \ref{characterization-theorem} (1)(2) in the general torsion free case:} 
\\
Suppose $(\nabla \gamma)(t_0)$ and $(\nabla^2\gamma)(t_0)$ are linearly independent. 
Then $f = \nabla{\mbox{\rm -}}\Tan(\gamma)$ is a frontal with the frame 
$V_1(t, s) = \frac{\pa f}{\pa t}(t, s)$ and $V_2(t, s) = F(t, s)$ for the integral lifting $\widetilde{f}$. 
We take coframe $L = (L_1, \dots, L_{m-2})$, $L_i : (I\times\R, (t_0, 0)) \to T^*M \setminus \zeta$ along $f$ satisfying 
$$
\langle L_i(t, s), V_1(t, s) \rangle = 0, \quad \langle L_i(t, s), V_2(t, s)\rangle = 0, \  (1 \leq i \leq m-2). 
$$
Set $\ell_i(t) = L_i(t, 0)$. 
Then we have 
$$
\langle \ell_i(t), (\nabla \gamma)(t) \rangle = 0, \quad \langle \ell_i(t), (\nabla^2\gamma)(t) \rangle = 0. 
$$
Set 
$$
\psi_i(t) := \langle \ell_i(t), (\nabla_\eta^f F)(t, 0) \rangle, (1 \leq i \leq m-2). 
$$
Since $\nabla$ is torsion free, by Lemma \ref{third derivative}, 
we have  $(\nabla_\eta^f F)(t, 0) = (\nabla^3\gamma)(t)$ 
and so 
$$
\psi_i(t) = \langle \ell_i(t), (\nabla^3\gamma)(t)\rangle. 
$$
Define the vector valued function 
$$
\psi : I \to \R^{m-2}, \quad \psi(t) = (\psi_1(t), \dots, \psi_{m-2}(t)). 
$$
Then $\psi(t_0) = 0$ if and only if $(\nabla\gamma)(t_0), (\nabla^2\gamma)(t_0), (\nabla^3\gamma)(t_0)$ 
are linearly dependent. 

Note that the covariant derivative $(\nabla_{\frac{\pa}{\pa t}}^\gamma \psi_i)(t)$ of the function $\psi_i(t)$ is equal to the 
ordinary derivative $\psi'_i(t)$. So we have 
$$
\psi_i'(t) = (\nabla_{\frac{\pa}{\pa t}}^\gamma \psi_i)(t) = \langle (\nabla_{\frac{\pa}{\pa t}}^\gamma \ell_i)(t), (\nabla^3\gamma)(t) \rangle 
+ \langle \ell_i(t), (\nabla^4 \gamma) (t)  \rangle. 
$$
Since $\langle \ell_i(t), (\nabla\gamma)(t)\rangle = \langle \ell_i(t), (\nabla^2\gamma)(t)\rangle = 0$, we have 
$$
0 = \langle (\nabla_{\frac{\pa}{\pa t}}^\gamma \ell_i)(t), (\nabla\gamma)(t) \rangle + \langle \ell_i(t), (\nabla^2\gamma)(t) \rangle = 
\langle (\nabla_{\frac{\pa}{\pa t}}^\gamma \ell_i)(t), (\nabla\gamma)(t) \rangle. 
$$

If $\psi(t_0) = 0$ namely if $(\nabla\gamma)(t_0), (\nabla^2\gamma)(t_0), (\nabla^3\gamma)(t_0)$ 
are linearly dependent, then we have that
$
\langle \ell_i(t_0), (\nabla^3\gamma)(t_0) \rangle = 0,  (1 \leq i \leq m-2), 
$
since $(\nabla^3\gamma)(t_0)$ is a linear combination of $(\nabla\gamma)(t_0)$ and $(\nabla^2\gamma)(t_0)$. 
Thus we have 
$$
0 = \langle(\nabla_{\frac{\pa}{\pa t}}^\gamma \ell_i)(t_0), (\nabla^2\gamma)(t_0) \rangle + \langle \ell_i(t_0), 
(\nabla^3\gamma)(t_0) \rangle = 
\langle(\nabla_{\frac{\pa}{\pa t}}^\gamma \ell_i)(t_0), (\nabla^2\gamma)(t_0) \rangle. 
$$
Moreover we have $\langle(\nabla_{\frac{\pa}{\pa t}}^\gamma \ell_i)(t_0), (\nabla^3\gamma)(t_0) \rangle = 0$. 
Therefore we have 
$$
\psi_i'(t_0) = \langle (\nabla_{\frac{\pa}{\pa t}}^\gamma \ell_i)(t_0), (\nabla^3\gamma)(t_0)\rangle + 
\langle \ell_i(t_0), (\nabla^4\gamma)(t_0) \rangle = \langle \ell_i(t_0), (\nabla^4\gamma)(t_0)\rangle. 
$$
Therefore, if $\psi(t_0) = 0$ and 
$(\nabla\gamma)(t_0), (\nabla^2\gamma)(t_0), (\nabla^4\gamma)(t_0)$ are linearly independent, 
then $\psi'(t_0) \not= 0$. 

Now, by Theorem \ref{characterization}, 
we see that, if $\dim(M) = m = 3$, and $\psi(t_0) = 0, \psi'(t_0) \not= 0$, 
then $f$ is diffeomorphic to the folded umbrella (cuspidal cross cap) 
at $(t_0, 0)$. 
Moreover, if $\dim(M) \geq 3$ and $\psi(t_0) \not= 0$, then $f$ is diffeomorphic to 
the cuspidal edge at $(t_0, 0)$. 
\QED

\

\section{Proof of genericity theorem}
\label{Proof of the genericity theorem 1}

In general, let $\gamma : I \to M$ be a $C^\infty$ curve and $t_0 \in I$. Define 
$$
a_1 := \inf\left\{ k \ \left\vert \ k \geq 1, (\nabla^k\gamma)(t_0) \not= 0 \right.\right\}. 
$$
Note that $\gamma$ is an immersion at $t_0$ if and only if $a_1 = 1$. 
If $a_1 < \infty$, then define
$$
a_2 := \inf\left\{ k \left\vert \ \rank\left( (\nabla\gamma)(t_0), (\nabla^2\gamma)(t_0), \dots, (\nabla^k\gamma)(t_0) \right) 
= 2 \right.\right\}. 
$$
We have $1 \leq a_1 < a_2$. If $a_i < \infty, 1 \leq i < \ell \leq m$, then define $a_\ell$ inductively by 
$$
a_\ell := \inf\left\{ k \ \left\vert \
\rank\left( (\nabla\gamma)(t_0), (\nabla^2\gamma)(t_0), \dots, (\nabla^k\gamma)(t_0) \right) 
= \ell \right.\right\}. 
$$
If $a_m < \infty$, then we call the strictly increasing sequence $(a_1, a_2, \dots, a_m)$ of 
natural numbers the {\it $\nabla$-type} of $\gamma$ at $t_0$.

To obtain our Theorem \ref{genericity-theorem}, we show

\bep
\label{genericity1}
Let $(M, \nabla)$ be a manifold of dimension $m$ with an affine connection $\nabla$. 
Then there exists an open dense set ${\mathcal U}$ in the set of $C^\infty$ curves from $I$ to $M$ 
in Whitney $C^\infty$ topology, such that for any $\gamma$ belonging to ${\mathcal U}$ and for any $t_0 \in I$, 
$\gamma$ is of $\nabla$-type 
$(1, 2, 3)$ or $(1, 2, 4)$ if $m = 3$, and 
$(1, 2, 3, \dots, m-1, m)$ or $(1, 2, 3, \dots, m-1, m+1)$ if $m \geq 4$, at $t_0$. 
\enp

\noindent
{\it Proof of Propositions \ref{genericity1}:} 
First we remark that, for any local coordinates on $M$, and 
for each $k = 1, 2, \dots$, the iterated covariant derivative
$(\nabla^{k}\gamma)(t)$ is expressed as 
$(\nabla^{k}\gamma)(t) = \gamma^{(k)}(t) + P$, by a polynomial $P$ of $\gamma^{(i)}(t), 0 \leq i < k$ and 
$(\pa^\alpha \Gamma^\lambda_{\mu\nu} / \pa x^\alpha)(\gamma(t)), \vert\alpha\vert \leq k - 2$ 
(cf. Lemma \ref{iteration}). Therefore, for positive integer $r$, 
there exists an algebraic diffeomorphism of $\Phi : J^r(I, M)_{t_0, q} \to J^r(I, M)_{t_0, q}$ 
of the $r$-jet space $J^r(I, M)_{t_0, q} = \{ j^r\gamma(t_0) \mid \gamma : (I, t_0) \to (M, q)\}$ 
satisfying the following conditions: 
if $\Phi(j^r\gamma(t_0)) = j^r\beta(t_0), \beta : (I, t_0) \to (M, q)$, then 
$(\nabla^k\gamma)(t_0) = \beta^{(k)}(t_0), 1 \leq k \leq r$. In particular, 
for any $1 \leq \ell \leq r$, and for any $1 \leq i_1 < i_2 < \cdots < i_\ell \leq r$ and for any $j, 0 \leq j \leq m$, 
$\rank((\nabla^{i_1}\gamma)(t_0), \dots, (\nabla^{i_\ell}\gamma)(t_0)) = j$ if and only if 
$\rank(\beta^{(i_1)}(t_0), \dots, \beta^{(i_\ell)}(t_0)) = j$. 
Let $r \geq m+1$. 
We set 
\begin{align*}
S_{\nabla} := \{ j^{r}\gamma(t_0) \ \mid \ & (\nabla \gamma)(t_0), (\nabla^2 \gamma)(t_0), 
\dots, (\nabla^{m}\gamma)(t_0) 
{\mbox{\rm \ are linearly dependent and }} 
\\
& 
(\nabla \gamma)(t_0), (\nabla^2 \gamma)(t_0), \dots, (\nabla^{m-1}\gamma)(t_0), (\nabla^{m+1}\gamma)(t_0) 
{\mbox{\rm \ are linearly dependent.}}\}
\end{align*}
is an algebraic set of codimension $2$. 
The above diffeomorphism $\Phi$ maps $S_{\nabla}$ to $S_{\nabla_0}$. Here $\nabla_0$ is the trivial connection on $\R^m$. 
Thus the calculation is reduced to the trivial case which is well known (see \cite{GG}\cite{Ishikawa4}). 
Note that  $S_{\nabla}$ is intrinsically defined by the given connection $\nabla$. 
Then we have the associated closed stratified subbundle
$S_{\nabla}(I, M)$ of the $r$-jet bundle $J^r(I, M) \to I\times M$ of codimension $2$. 
By the transversality theorem 
$$
{\mathcal U} := \{ \gamma : I \to M \mid j^r\gamma : I \to J^r(I, M) {\mbox{\rm \ is transverse to }} S_{\nabla}(I, M) 
\} 
$$
is open dense in Whitney $C^\infty$ topology (see \cite{GG}). 
Let $\gamma \in {\mathcal U}$ and $t_0 \in I$. Since $S_{\nabla}(I, M)$ is codimension $2$, 
$j^r\gamma(t_0) \not\in S_{\nabla}(I, M)$. This means that 
$$
(\nabla \gamma)(t_0), (\nabla^2 \gamma)(t_0), \dots, (\nabla^{m}\gamma)(t_0)
$$
are linearly independent, or, they are linearly dependent but 
$$
(\nabla \gamma)(t_0), (\nabla^2 \gamma)(t_0), \dots, (\nabla^{m-1}\gamma)(t_0), (\nabla^{m+1}\gamma)(t_0)
$$
are linearly independent. In the first case, $\gamma$ is of $\nabla$-type $(1, 2, \dots, m-1, m)$. 
In the second case, $\gamma$ is of type $\nabla$-type $(1, 2, \dots, m-1, m+1)$. 
Thus we have Proposition \ref{genericity1}. 
\QED

\ber
{\rm 
Let $\aaa = (a_1, a_2, \dots, a_m)$ be any strictly increasing sequence of positive integers. 
For an integer $r \geq a_m$, we define, in the $r$-jet bundle $J^r(I, M)$, 
$$
\Sigma_{\aaa}(I, M) : = \{ j^r\gamma(t_0) \in J^r(I, M) \mid 
\gamma {\mbox{\rm \ is of \ }} \nabla{\mbox{\rm -type\ }} \aaa \}. 
$$
Then $\Sigma_{\aaa}(I, M)$ is a stratified subbundle of $J^r(I, M)$ over $I\times M$ with 
an algebraic typical fiber. It can be shown, as in the proof of Propositions \ref{genericity1}, 
that the codimension of 
$\Sigma_{\aaa}(I, M)$ in $J^r(I, M)$ is independent of $\nabla$ and 
is given by $\sum_{i = 1}^m (a_i - i)$ (see \cite{Ishikawa4} for the flat case). 
}
\enr

\noindent
{\it Proof of Theorem \ref{genericity-theorem}:}
We may suppose $\nabla$ is torsion free (Remark \ref{torsion-free}). Then Proposition \ref{genericity1}, 
and Theorem \ref{characterization-theorem} imply Theorem \ref{genericity-theorem}.

\section{Tangent surfaces to torsionless curves and fold singularities}
\label{Tangent surfaces to torsionless curves and fold singularities}

Let $(M, \nabla)$ be a manifold of dimension $m \geq 3$ 
with a torsion free affine connection $\nabla$.  
Consider a curve $\gamma : I \to M$ such that $(\nabla\gamma)(t), (\nabla^2\gamma)(t)$ are linearly 
independent for any $t \in I$. Though the torsion is not defined in general, we can define: 

\bef
{\rm 
A curve $\gamma : I \to M$ is called {\it torsionless} if 
$(\nabla\gamma)(t), (\nabla^2\gamma)(t)$ are linearly 
independent but 
$(\nabla\gamma)(t), (\nabla^2\gamma)(t), (\nabla^3\gamma)(t)$ 
are linearly dependent everywhere. 
}
\enf

The situation is, by any means, non-generic. However 
to study torsionless curves is an interesting geometric problem (cf. \cite{Cartan}). 

\ber
{\rm 
Let $M$ be a Riemannian manifold with the Levi-Civita connection $\nabla$. The notion of torsion is well-defined for the curve $\gamma$ parametrized by the arc-length, if $\nabla\gamma, \nabla^2\gamma$ are linearly independent. 
Then the condition that the torsion of $\gamma$ is zero if and only if $\gamma$ satisfies the third order non-linear ordinary differential equation, 
$$
\nabla^3\gamma - \frac{\Vert\nabla^2\gamma\Vert'}{\Vert\nabla^2\gamma\Vert} \, \nabla^2\gamma
+ \Vert\nabla^2\gamma\Vert^2 \, \nabla\gamma = 0. 
$$
Then in particular $\gamma$ is torsionless in our sense. 
}
\enr

If $(M, \nabla)$ is projectively flat, namely, if it is 
projectively equivalent to the Euclidean space $(E^m, \nabla_0)$ 
with the standard connection $\nabla_0$, then it is well-known that any torsionless curve is a \lq\lq plane curve", 
therefore its $\nabla$-tangent surface is \lq\lq folded" into a totally geodesic surface. 
Note that a Riemannian manifold is locally projectively flat if and only if it has a constant curvature 
(Beltrami's Theorem, see p.352 of \cite{Sharpe} or p.97 of \cite{Eisenhart}). 

A map-germ $f : (\R^2, p) \to (M, f(p))$ is called a {\it fold}, or an {\it embedded fold}, if it is diffeomorphic to 
$$
(t, s) \mapsto (t + s, t^2 + 2st, 0, \dots, 0), 
$$
which is diffeomorphic also to $(u, w) \mapsto (u, \frac{1}{2}w^2, 0, \dots, 0)$. 
The fold singularities appear in other geometric problems also (see \cite{FRUYY} for instance). 
For our problem, we have 

\bep
\label{tangent-to-torsionless-curve}
Let $(M, \nabla)$ be locally projectively flat around $q \in M$ and  
$\gamma : (\R, t_0) \to (M, q)$ a germ of torsionless curve. 
Then 
the germ of $\nabla$-tangent surface $\nabla{\mbox{\rm -}}\Tan(\gamma) : (\R^2, (t_0, 0)) \to 
(M, q)$ to $\gamma$ is a fold. 
In particular it is a generically two-to-one mapping. 
\enp

\beP
There exits a germ of projective equivalence $\varphi : (M, q) \to (E^m, 0)$ such that 
$\varphi\circ \gamma : (\R, 0) \to (E^m, 0)$ is a strictly convex plane curve in $E^2\times \{ 0\} \subset E^m$. 
Then, by Theorem \ref{characterization2}, 
the tangent surface 
$\nabla{\mbox{\rm -}}\Tan(\varphi\circ\gamma) = \varphi\circ(\nabla{\mbox{\rm -}}\Tan(\gamma))$ 
is diffeomorphic to a fold regarded as a map-germ $(\R^2, (0, 0)) \to (\R^2, (0, 0))$. 
Therefore $\nabla{\mbox{\rm -}}\Tan(\gamma)$ is diffeomorphic to a fold, 
regarded as a surface-germ $(\R^2, (t_0, 0)) \to (M, q)$. 
\QED

\

A map-germ $f : (\R^2, p) \to (M^3, f(p))$ to a three dimensional space 
is called a {\it $(2, 5)$-cuspidal edge}, if it is diffeomorphic to 
$$
(u, w) \mapsto (u, w^2, w^5). 
$$
Then $f$ is a frontal with a non-degenerate singular point at $p$,  
and the characteristic function $\psi$ vanishes identically (see \S \ref{Cuspidal edge and folded umbrella}). 
In fact, the model germ
$(u, w) \mapsto (u, w^2, w^5)$ 
(the \lq\lq suspension" of $(2, 5)$-cusp) is a frontal 
with a non-degenerate singular point at $0$, the kernel field $\pa/\pa u$ for $f_*$ is transverse to the singular curve $w = 0$, and 
the characteristic function $\psi(u) \equiv 0$. 
However $f$ is injective and it is never diffeomorphic to a tangent surface to any 
torsionless curve in a projectively flat space, by Proposition \ref{tangent-to-torsionless-curve}. 
Nevertheless it can be a $\nabla$-tangent surface of a torsionless curve 
for a torsion free affine connection on $M^3$. 

\bee
{\rm 
Let $\nabla$ be the torsion free affine connection on $\R^3$ with coordinates $x_1, x_2, x_3$ 
defined by 
$\Gamma^\lambda_{\mu\nu} = x_1 + x_2^2$ if $(\lambda, \mu, \nu) = (3, 1, 2), (3, 2, 1)$ and 
otherwise $\Gamma^\lambda_{\mu\nu} = 0$. Let $\gamma : \R \to \R^3$ 
be the immersion defined by $\gamma(t) = (-t^2, t, 0)$. Then $\Gamma^\lambda_{\mu\nu} = 0$ along 
$\gamma$. We have $(\nabla\gamma)(t) = (-2t, 1, 0), (\nabla^2\gamma)(t) = (-2, 0, 0), 
(\nabla^3\gamma)(t) = (0, 0, 0)$. Therefore $\gamma$ is torsionless. 
For any $t_0 \in \R$, the $\nabla$-geodesic with the initial condition $(\gamma(t_0), \gamma'(t_0))$ 
is given by $(-2t_0s - t_0^2, s + t_0, \frac{1}{3}t_0 s^4)$. Therefore we have 
$$
f(t, s) = \nabla{\mbox{\rm -}}\Tan(\gamma)(t, s) = (-2ts - t^2, \ s+t, \ \frac{1}{3}ts^4). 
$$
Set $u = s+t, w = s$. Then we see that $f$ is diffeomorphic to 
$(u, w) \mapsto (-u^2+w^2, u, \frac{1}{3}uw^4 - \frac{1}{3}w^5)$, which is diffeomorphic to 
$(u, w) \mapsto (u, w^2, w^5)$. Therefore $f$ is a $(2, 5)$-cuspidal edge. 
}
\ene

We conclude this section by posing the problem on singularities of $\nabla$-tangent surfaces 
to torsionless curves in the case of a general affine connection $\nabla$.

\

\begin{flushleft}
Goo ISHIKAWA, \\
Department of Mathematics, Hokkaido University, 
Sapporo 060-0810, Japan. \\
e-mail : ishikawa@math.sci.hokudai.ac.jp \\

\

Tatsuya YAMASHITA, \\
Department of Mathematics, Hokkaido University, 
Sapporo 060-0810, Japan. \\
e-mail : tatsuya-y@math.sci.hokudai.ac.jp \\ 

\end{flushleft}

\end{document}